\documentclass[a4paper,11pt,twoside]{article}

\usepackage{xypic}
\usepackage{latexsym}
\usepackage{amsmath}
\usepackage{amstext}
\usepackage{amsfonts}
\usepackage{amssymb}
\usepackage{amsbsy}
\usepackage{eucal}
\usepackage{theorem}
\xyoption{all}
\CompileMatrices

\newtheorem{theorem}{Theorem}
\newtheorem{lemma}{Lemma}
\newtheorem{proposition}[lemma]{Proposition}
\newtheorem{corollary}[lemma]{Corollary}

\theorembodyfont{\rmfamily}
\newtheorem{remark}[lemma]{Remark}

\newcommand{\Res}{\text{\rm Res}}

\newcommand{\qed} {\hbox{} \nolinebreak \hfill $\;\Box$}

\begin{document}

\addtolength{\baselineskip}{0.5 pt}

\setlength{\parskip}{0.2 ex}

\title{On the Structure of Weil Restrictions of \\ Abelian Varieties}

\author{Claus Diem \and Niko Naumann}

\date{June 10, 2003}

\maketitle

\pagestyle{myheadings}
\markboth{\sc Claus Diem, Niko Naumann}{\sc On the Structure of Weil Restrictions of Abelian Varieties}

\begin{abstract}
We give a description of endomorphism rings of Weil restrictions of abelian varieties with respect to finite Galois extensions of fields. The results are applied to study the isogeny decompositions of Weil restrictions.
\end{abstract}


\section*{Introduction}
For the use of Weil restrictions of abelian varieties in various fields of mathematics but also because of genuine interest in Weil restrictions themselves, it is important to determine the endomorphism rings and the isogeny decompositions. This is what this article provides -- at least in two important special cases.

After giving a brief expos\'e of general facts about Weil restrictions of abelian varieties 
in the first section, we study Weil restric\-tions with respect to extensions of finite fields in the second section. Here we determine the endomorphism algebra of a Weil restriction (see Theorem \ref{finite-end-structure}) and then show that under rather general assumptions, the Weil restriction is simple over the base-field (see Theorem \ref{new-simple-res-simple}).

In the third section, we deal with the following situation: $K|k$ is an arbitrary finite Galois extension of fields, $A$ an abelian variety over $k$, $W$ the Weil restriction of $A_K$ with respect to $K|k$. We describe the endomorphism ring of $W$ as a skew group ring over $\text{End}(A_K)$ (see Theorem \ref{theorem-skew-group-ring-end-W-iso}) and apply this result to study the isogeny decomposition of $W$ over $k$. In the last subsection the results are applied to give an explicit description of the isogeny decomposition of $W$ in the case of a cyclic field extension; see Theorem \ref{W-factors}.

\subsection*{Notation}

\subsubsection*{General}
By a \emph{ring} we mean a ring with unity, and by a ring-homomorphism a homomorphism of rings with unity. If $R$ is a ring and $\Sigma$ a finite set, then by $\text{\rm M}_\Sigma(R)$ we mean the matrix ring over $R$ on the set $\Sigma$. For any abelian group $G$, $G^{\circ}$ denotes $G\otimes\mathbb{Q}$. If $D$ is a skew field, we denote its center by $Z(D)$. 
If $k$ is a field, $\overline{k}$ denotes an algebraic closure. If $k$ is some field and $X$ and $Z$ are $k$-schemes, we denote the $k$-morphisms from $Z$ to $X$ by $X(Z)$.

Let $k$ be a field. By a \emph{homomorphism} between abelian $k$-varieties we mean a  morphism of $k$-schemes which preserves the group structure (other authors might call this a $k$-homomorphism or a $k$-morphism of abelian varieties). Analogous definitions apply to isogenies and endomorphisms. The group of homomorphisms between two abelian $k$-varieties $A$ and $B$ is denoted by $\text{Hom}(A,B)$ and the ring of endomorphisms of an abelian $k$-variety $A$ by $\text{End}(A)$. Following this terminology, we use the notion of a \emph{simple} abelian $k$-variety where other authors might speak of a $k$-simple abelian $k$-variety. If two abelian $k$-varieties $A$ and $B$ are isogenous, we write $A \sim B$.

If we are given an extension of fields $K|k$, we denote $k$-schemes by $X,Y$ etc.\ and $K$-schemes by $X',Y'$ etc.\ (or by $X_K,Y_K$ etc.\ if they are induced by base-change $K|k$).

We denote the dual abelian variety of an abelian $k$-variety $A$ by $\widehat{A}$. For an invertible sheaf ${\cal L}$ on $A$, $\phi_{\cal L}: A \longrightarrow \widehat{A}$ denotes the corresponding homomorphism; c.f. \cite[\S 6]{Mu}. Following \cite{Mi-2}, a polarization $\varphi$ of $A$ is a homomorphism $A \longrightarrow \widehat{A}$ such that $\varphi \otimes_k \text{id}_{\overline{k}} = \phi_{\mathcal{L}} : A_{\overline{k}} \longrightarrow \widehat{A}_{\overline{k}}$ for some ample invertible sheaf on $A_{\overline{k}}$.

\subsubsection*{Galois twists}
Let $K|k$ be a Galois extension of fields with Galois group $G$. Then the elements of $G$ induce automorphisms of the $\text{Spec}(k)$-scheme $\text{Spec}(K)$ -- we obtain in this way an anti-isomorphism $G \longrightarrow \text{Aut}_{\text{Spec}(k)}(\text{Spec}(K))$.

We identify the opposite group $G^{\text{opp}}$ with $\text{Aut}_{\text{Spec}(k)}(\text{Spec}(K))$. We will always work with $G^{\text{opp}}$ instead of $G$.

Let $X'$ be a $K$-scheme.

For $\sigma \in G^{\text{opp}}$, let $\sigma^{-1}(X')$ be the pull-back of $X'$ via $\sigma : \text{Spec}(K) \longrightarrow \text{Spec}(K)$, i.e.\ if $p_K : X' \longrightarrow \text{Spec}(K)$ is the structure morphism, $\sigma^{-1}({X'})$ is $X'$ considered as $K$-scheme via $\sigma^{-1} \circ p_K$. We denote the canonical isomorphism of $k$-schemes from $\sigma^{-1}(X')$ to $X'$ also by $\sigma$. If $Y'$ is another $K$-scheme and $\alpha: X' \longrightarrow Y'$ is a morphism of $K$-schemes, we obtain by base-change a morphism of $K$-schemes $\sigma^{-1}(\alpha) = \sigma^{-1} \alpha \sigma : \sigma^{-1}(X') \longrightarrow \sigma^{-1}({Y'})$.

If $X'$ is an abelian $K$-variety, by pull-back $\sigma^{-1}(X')$ also has the structure of an abelian $K$-variety.

\subsubsection*{Frobenius morphisms}
Let $q$ be a power of a prime number, $k$ the finite field with $q$ elements, let $A$ be an abelian $k$-variety. The \emph{Frobenius endomorphism} $\pi_k$ of $A$ is defined by the identity on the underlying topological space and by $f \mapsto f^q$ on the structure-sheaf $\mathcal{O}_{A}$. As the name indicates, $\pi_k$ is an endomorphism of the abelian $k$-variety $A$.

Now let $K|k$ be an algebraic extension of fields. We identify the Galois group Gal$(K|k)$ with its dual. The Frobenius automorphism of $K|k$ (or of $\text{Spec}(K) \longrightarrow \text{Spec}(k)$) is denoted by $\sigma_{K|k}$. If $K = \overline{k}$, we write $\sigma_k$ instead of $\sigma_{\overline{k}|k}$.

Let $A'$ be an abelian $K$-variety. As stated above, we have a canonical isomorphism of $k$-schemes $\sigma_{K|k} : \sigma_{K|k}^{-1}(A') \longrightarrow A'$. The \emph{relative Frobenius homomorphism} (with respect to $k$) $\pi_k : A' \longrightarrow \sigma_{K|k}^{-1}(A')$ is a homomorphism of abelian $K$-varieties which is defined as follows: Let $F_{k}$ be the morphism of the $k$-scheme $A'$ to itself which is the identity on the underlying topological space and it is given by $f \mapsto f^q$ on the structure-sheaf $\mathcal{O}_{A'}$. Then $\pi_k := \sigma_{K|k}^{-1} \circ F_{k} : A' \longrightarrow \sigma_{K|k}^{-1}(A')$.

\section{Definitions and first results}

\subsection{Definition of the Weil restriction}
Let $K|k$ be a finite Galois extension. Let  $A'$ an abelian $K$-variety. It is well-known that the functor
\[ Z \mapsto A'(Z \otimes_k K)
\]
from the category of $k$-schemes to the category of abelian groups is representable by an abelian $k$-variety; for a construction via Galois theory see Subsection \ref{construction}, for a construction via ``restriction of scalars'' see \cite[7.6]{BLR}. (The representatility of the functor by an abelian variety holds more generally for a finite separable extension of fields, but we restrict ourselves to the Galois-case is this article.) A representing object will be denoted $\Res_k^K(A')$ and will be called the \emph{Weil restriction} of $A'$ with respect to $K|k$. The universal element $u \in A'(\Res_k^K(A') \otimes_k K)$ maps the zero of $\Res_k^K(A') \otimes_k K$ to the zero of $A'$ and thus is a homomorphism of abelian $K$-varieties.

Now, $\Res_k^K(A')$ with $u$ is also a representing object for the functor $B \mapsto \text{Hom}(B_K,A')$ from the category of abelian $k$-varieties to the category of abelian groups as well as for the functor $B \mapsto \text{Hom}^{\! \circ \!}(B_K,A')$ from the category of abelian $k$-varieties up to isogeny to the category of $\mathbb{Q}$-vector spaces.

\subsection{Construction of the Weil restriction}
\label{construction}
Let us recall the construction of $\Res_k^K(A')$ via Galois theory.

Let $W'$ be the following product of Galois-conjugates of $A'$:
\begin{equation}
\label{Product-galois-conjugates}
W' := \prod_{\sigma \in G^{\text{opp}}} \sigma^{-1}(A')
\end{equation}
Let $p_\sigma : W' \longrightarrow \sigma^{-1}(A')$ be the projections, let $\text{Aut}_k(W')$ be the group of automorphisms of the $k$-scheme $W'$.

We define a Galois operation on $W'$ by $G^{\text{opp}} \longrightarrow \text{Aut}_k(W'), \; \tau \mapsto \widetilde{\tau}$ where $\widetilde{\tau} = (\tau \, p_{\sigma \tau})_{\sigma \in G^{\text{opp}}}$.
Since $W'$ is projective, the quotient $W:=W'/G$ under this operation exists and is projective. We have $W' \simeq W_K$.

Fix some $k$-scheme $Z$. We have a Galois operation on $W'(Z  \otimes_k  K)$. If $\tau \in G^{\text{opp}}$ and $P=(P_{\sigma})_{\sigma \in G^{\text{opp}}} \in W'(Z \otimes_k K)$, then $\tau((P_{\sigma})_{\sigma \in G^{\text{opp}}}) = (\tau(P_{\sigma \tau}))_{\sigma \in G^{\text{opp}}}$.
It follows that $P \mapsto (\sigma^{-1}(P))_{\sigma \in {G^{\text{opp}}}}$ is a bijection between the $Z \otimes_k K$-valued points of $A'$ and the Galois-invariant $Z \otimes_k K$-valued points of $W'$. On the other hand, by Galois theory, the Galois-invariant $Z \otimes_k K$-valued points of $W'$ are in bijection with the $Z$-valued points of $W$. Both bijections are natural in $Z$.

It follows that $W = W'/G$ with universal element $u := p_{\text{id}}$ represents the functor $Z \mapsto A'(Z \otimes_k K)$ from the category of $k$-schemes to the category of sets. Via the group laws on these sets, one defines a group law on $W$, and with this group law, $W$ is an abelian variety. By construction, the neutral element and the addition law of $W$ coincide after base-change with the neutral element and the addition law of the product of Galois-conjugates in (\ref{Product-galois-conjugates}). Moreover, the universal element $u = p_{\text{id}}$ is a homomorphism of abelian $K$-varieties.

From the Galois-operation of $W'$, we obtain 
\begin{equation}
\label{operation-on-projection}
\tau^{-1}(p_\sigma) = p_{\sigma \tau}, \text{ especially }
\tau^{-1}(u) = p_\tau.
\end{equation}

\subsection{The functor ``restriction of scalars''}
\label{functor-Res}
The assignment $A' \mapsto \Res_k^K(A')$ defines a covariant additive functor $\text{Res}_k^K$ from the category of abelian $K$-varieties (up to isogeny) to the category of abelian $k$-varieties (up to isogeny). This functor is called ``restriction of the field of definition'' or ``restriction of scalars'' or ``norm functor''; cf.\ \cite{Mi-AA}.

For any abelian $K$-variety $A'$, $\Res^K_k$ gives a ring-homomorphism from $\text{\rm End}(A')$ to $\text{\rm End}(\Res^K_k(A'))$ and from $\text{\rm End}^{\! \circ \!}(A')$ to $\text{\rm End}^{\! \circ \!}(\Res^K_k(A'))$.

Let $A', B'$ be abelian $K$-varieties. Then
\begin{equation} \label{Hom-Res}
\text{\rm Hom}(\Res^K_k(A')_K,\Res^K_k(B')_K) 
\simeq \bigoplus_{\sigma, \nu \in {G^{\text{opp}}}} \text{\rm Hom}(\nu^{-1}({A'}), \sigma^{-1}({B'}));
\end{equation}
see equations (\ref{Product-galois-conjugates}) and (\ref{prod-morphisms-matrix}). 
\\

Let $\alpha : A' \longrightarrow B'$ be a homomorphism. Then under (\ref{Hom-Res}), $\Res^K_k(\alpha)\otimes_k \text{id}_K$ is given by the diagonal ``matrix'' 
\[(\sigma^{-1}(\alpha) \delta_{\sigma,\nu})_{\sigma, \nu \in {G^{\text{opp}}}} \in \bigoplus_{\sigma, \nu \in {G^{\text{opp}}}} \text{\rm Hom}(\nu^{-1}({A'}), \sigma^{-1}({B'})),\]
where $\delta_{\sigma,\nu}$ is the ``Kronecker delta''. If $\alpha : A' \longrightarrow B'$ is an isogeny, then $\Res^K_k(\alpha) : \Res^K_k(A') \longrightarrow \Res^K_k(B')$ is an isogeny of degree $(\text{deg}(\alpha))^{[K:k]}$.

\subsection{The Weil restriction of the dual abelian variety}
The Weil restriction of the dual abelian variety is functorially isomorphic to the dual abelian variety of the Weil restriction. This can be seen as follows.

Let $W:= \Res_k^K(A')$.

Let $Z$ be some $k$-scheme, $\mathcal{L}$ some invertible sheaf on $A' \times_K Z_K$, algebraically equivalent to zero. Now consider the invertible sheaf
\[ \mathcal{L}_{W_K} := \bigotimes_\sigma p_\sigma^* \sigma^*(\mathcal{L}) = \bigotimes_\sigma \widetilde{\sigma}^* u^*(\mathcal{L})\]
on $W_K$. The isomorphism class in $\text{Pic}(W_K \times_K Z_K)/\text{Pic}(Z_K)$ of this invertible sheaf corresponds to an element in $\widehat{W}_K(Z_K)$ which is invariant under the Galois-operation and thus defines an element in $\widehat{W}(Z)$.

We obtain in this way a homomorphism $\widehat{A'}(Z_K) \longrightarrow \widehat{\Res_k^K(A')}(Z)$ which is functorial in $Z$. We thus have a homomorphism $\Res(\widehat{A'}) \longrightarrow \widehat{\Res^K_k(A')}$.

After base-change $K|k$, this homomorphism becomes the canonical isomorphism
\[ \prod_{\sigma \in G^{\text{opp}}} \sigma^{-1}(\widehat{A'}) \longrightarrow \widehat{\prod_{\sigma \in G^{\text{opp}}} \sigma^{-1}(A')}, \]
thus it is an isomorphism. This isomorphism $\Res(\widehat{A'}) \longrightarrow \widehat{\Res^K_k(A')}$ is functorial in $A'$ as can for example easily be seen after base-change $K|k$. We thus have:
\begin{proposition}
For abelian $K$-varieties $A'$, $\Res_k^K(\widehat{A'})$ is functorially isomorphic to $\widehat{\Res_k^K(A')}$.
\end{proposition}

\subsection{Weil restrictions of polarized abelian varieties}
\label{Weil restriction of pol. abelian varieties}
Let $K|k$ be a finite Galois field extension, $A'$ an abelian $K$-variety, $\widehat{A'}$ the dual abelian variety.

Let $\varphi: A' \longrightarrow \widehat{A'}$ be a polarization of $A'$, defined by an ample invertible sheaf $\mathcal{L}$ on $A'_{\overline{K}}$, i.e.\ $\varphi \otimes_K \text{id}_{\overline{K}}= \phi_{\mathcal{L}} : A'_{\overline{K}} \longrightarrow \widehat{A'}_{\overline{K}}$.
As stated in Subsection \ref{functor-Res}, this induces an isogeny
\[\Res^K_k(\varphi) : \, \Res^K_k(A') \longrightarrow \Res^K_k(\widehat{A'}) \simeq \widehat{\Res^K_k(A')}. \]
We show now that this homomorphism is again a polarization.

Let $\sigma \in G^{\text{opp}}$. We regard $\sigma^{-1}(\widehat{A'})$ as the dual abelian variety of $\sigma^{-1}(A')$. 

Let ${\sigma'}$ be a $\text{Spec}(\overline{K})$-automorphism with $\pi\circ\sigma'=\sigma$ for the natural map $\pi:\text{Spec}(\overline{K})\rightarrow \text{Spec}(K)$. Then
\[\sigma^{-1}(\varphi) \otimes_K \text{id}_{\overline{K}} = {\sigma'}^{-1}(\phi_{\mathcal{L}}) = \phi_{{\sigma'}^*(\mathcal{L})}.\]
Here, the first equation is obvious by the definition of $\sigma'$ and the second equation is a general fact for all polarizations on abelian varieties. It can be checked rather easily on $\overline{K}$-valued points.

After base-change, we get
\[\Res^K_k(\varphi) \otimes_k \text{id}_K = (\sigma^{-1}(\varphi) \circ p_\sigma)_{\sigma \in {G^{\text{opp}}}} : \prod_{\sigma \in {G^{\text{opp}}}} \sigma^{-1}(A') \longrightarrow \prod_{\sigma \in {G^{\text{opp}}}} \sigma^{-1}(\widehat{A'}).\]
This is a product polarization defined by the ample invertible sheaf
\begin{equation}
\label{product-sheaf}
\mathcal{L}_{W_{\overline{k}}} := \bigotimes_\sigma (p_\sigma \otimes_K \text{id}_{\overline{K}})^* {\sigma'}^*(\mathcal{L})
\end{equation}
on $W_{\overline{k}}$. 

If one starts with an ample invertible sheaf $\mathcal{L}$ on $A'$, then analogously to (\ref{product-sheaf}), one defines an ample invertible sheaf $\mathcal{L}_{W_{K}}$ on $W_K$. The class of this sheaf in the Picard group is invariant under the operation of $\text{Gal}(K|k)$ and thus defines an ample invertible sheaf on $W$ (because the Picard functor of an abelian variety is representable) -- alternatively, one can also define explicitely a descent-datum on $\mathcal{L}_{W_K}$.

\begin{proposition}
\label{Res.Pol.Functor}
Let $K|k$ be a finite Galois field extension, $A'$ an abelian $K$-variety. If $\varphi$ is a polarization on $A'$ (defined by a sheaf on $A'$), then $\Res^K_k(\varphi)$ is a polarization on $\Res_k^K(A')$ (defined by a sheaf on $\Res_k^K(A')$). Furthermore $\text{\rm deg}(\Res^K_k(\varphi))=(\text{\rm deg}(\varphi))^{[K:k]}$.

Thus ``restriction of scalars'' is a functor from the category of polarized abelian $K$-varieties (with polarizations defined by sheaves on $A'$)
 to the category of polarized abelian $k$-varieties (with polarizations defined by sheaves on $\Res_k^K(A')$) which preserves principal polarizations.
\end{proposition}

\subsection{Appendix to Section 1: Products and the Rosati involution}
\label{Products-Rosati}
Let $k$ be a field, let $B_i$ for $i=1, \ldots, m$ and $A_j$ for $j=1, \ldots, n$ be abelian $k$-varieties. Let $A := \prod_{j=1, \ldots, n} A_j, \; B := \prod_{i=1, \ldots, m} B_i$. Let $\iota^A_j : A_j \longrightarrow A$ be the inclusions and let $p^A_j : A \longrightarrow A_j$ be the projections. (Similar definitions for $B$ as well as the corresponding dual abelian varieties $\widehat{A}$ and $\widehat{B}$.) Then 
\begin{equation} 
\label{prod-morphisms-matrix}
\begin{array}{ccc}
\text{\rm Hom}(A,B) & \longrightarrow & \bigoplus_{i,j} \text{\rm Hom}(A_j,B_i) \\
\psi & \mapsto & (p^B_i \psi \iota^A_j)_{i=1, \ldots, m, \, j=1, \ldots, n}
\end{array} \end{equation}
is an isomorphism. (The same is true for the corresponding groups \linebreak $\text{\rm Hom}^{\! \circ \!}(\ldots, \ldots)$ of both sides.)

Thus every homomorphism from $A$ to $B$ is uniquely determined by its ``matrix'', and conversely, every ``matrix'' determines a homomorphism. Furthermore, the composition of homomorphisms corresponds to the usual multiplication of matrices.

In particular, under (\ref{prod-morphisms-matrix}), $\text{\rm End}(A)$ is isomorphic to the ``matrix ring'' \linebreak $\bigoplus_{i,j} \text{\rm Hom}(A_j,A_i)$.

For later use we want to study how the Rosati involution with respect to a product polarization operates on the ``matrices''. It is convenient to generalize the concept of a ``Rosati involution'' first.

Let $X$ and $Y$ be abelian $k$-varieties with fixed polarizations $\varphi_X : X \longrightarrow \widehat{X}, \; \varphi_Y : Y \longrightarrow \widehat{Y}$. Then for every $\psi \in \text{\rm Hom}^{\! \circ \!}(X, Y)$, we denote $\varphi_X^{-1} \, \widehat{\psi} \, \varphi_Y \in \text{\rm Hom}^{\! \circ \!}(Y,X)$ by $\psi'$ and call it the \emph{Rosati involution} of $\psi$ with respect to $\varphi_X$ and $\varphi_Y$. 

Now for $i=1, \ldots, m, \; j=1, \ldots, n$, let $\varphi_{B_i} : B_i \longrightarrow \widehat{B_i}$ and $\varphi_{A_j} : A_i \longrightarrow \widehat{A_j}$ be polarizations. Let $\varphi_A$ and $\varphi_B$ be the corresponding product polarizations.
\begin{lemma}
\label{matrix-Rosati}
Let $\psi \in \text{\rm Hom}^{\! \circ \!}(A, B)$ be given by the ``matrix'' $(\psi_{i,j})_{i=1, \ldots, m, \; j=1, \ldots, n},$\linebreak$\psi_{i,j} \in \text{\rm Hom}^{\! \circ \!}(A_j,B_i)$. Then with respect to $\varphi_A$ and $\varphi_B$, the Rosati involution of $\psi$ is given by the ``matrix'' $(\psi_{j,i}')_{i=1, \ldots, n, \; j=1, \ldots, m}$ with  $\psi_{j,i}' \in \text{\rm Hom}^{\! \circ \!}(B_j,A_i)$.
\end{lemma}
\emph{Proof} Straightforward calculation. \qed

\section{Results for finite fields}
\label{new-factors}
Let $K|k$ be a finite extension of \emph{finite} fields of degree $n$. Let $A'$ be an abelian variety over $K$, 
$W$ the Weil restriction of $A'$ with respect to $K|k$.

\subsection{The endomorphism algebra}
We now study the endomorphism algebra and the isogeny decomposition of $W$ over $k$.

Let $\pi_k :  A' \longrightarrow \sigma_{K|k}^{-1}(A')$ be the relative Frobenius homomorphism with respect to $k$ and let $\pi_k : W \longrightarrow W$ be the Frobenius endomorphism; cf.\;``Notation''.

Let $\pi_K$ be the Frobenius endomorphism of $A'$. Then the image of $\pi_K$ under the ring-homomorphism $\Res^K_k$ equals the endomorphism $\pi_k^n$ of $W$. (In fact, after base-change $K|k$, $\Res^K_k(\pi_K)$ as well as $\pi_k^n$ become equal to the Frobenius endomorphism of $W_K$.) Thus the ring-homomorphism $\Res^K_k : \text{End}(A') \longrightarrow \text{End}(W)$ restricts to an inclusion $\mathbb{Z}[\pi_K] \longrightarrow \text{End}(W)$, given by $\pi_K \mapsto \pi_k^n$. This ring-homomorphism extends to a ring-homomorphism $\mathbb{Z}[\pi_K][X]/(X^n-\pi_k) \longrightarrow \text{End}(W)$, given by $X \longrightarrow \pi_k$.

The Frobenius endomorphism $\pi_k$ of $W$ commutes with all endomorphisms of $W$. Thus by the universal property of the tensor product, the ring-homomorphisms $\text{\rm End}(A') \longrightarrow \text{\rm End}(W), \; \lambda \mapsto \Res^K_k(\lambda)$ and $\mathbb{Z}[\pi_K][X]/(X^n-\pi_K) \longrightarrow \text{\rm End}(W), \; X \mapsto \pi_k$ induce a ring-homomorphism
\[ \text{\rm End}(A') \otimes_{\mathbb{Z}[\pi_K]} \mathbb{Z}[\pi_K][X]/(X^n-\pi_K) \longrightarrow \text{\rm End}(W), \; \lambda \mapsto \Res^K_k(\lambda) , \; X \mapsto \pi_k. \]

\begin{theorem}
\label{finite-end-structure}
Let $K|k$ be an extension of degree $n$ of finite fields. Let $A'$ be an abelian $K$-variety, $W$ the Weil restriction of $A'$ with respect to $K|k$. Then
\[
\text{\rm End}^{\! \circ \!}(A') \otimes_{\mathbb{Q}[\pi_K]} \mathbb{Q}[\pi_K][X]/(X^n-\pi_K) \longrightarrow \text{\rm End}^{\! \circ \!}(W), \; \lambda \mapsto \Res^K_k(\lambda) , \; X \mapsto \pi_k \]
is an isomorphism.
\end{theorem}
\emph{Proof}
By the defining property of the Weil restriction, as abelian groups, 
\begin{equation}
\label{finite-W-valued-points}
\text{\rm Hom}^{\! \circ \!}(W,W) \simeq \text{\rm Hom}^{\! \circ \!}(\prod_{i=0}^{n-1} \sigma_{K|k}^{-i}(A'),A') \text{ via } a \mapsto p_{\text{id}} \circ (a \otimes_k \text{id}_K).
\end{equation}
We show that the homomorphism of abelian groups
\begin{equation}
\begin{array}{c}
\label{finite-tensor-iso-row-vector}
\text{\rm Hom}^{\! \circ \!}(A',A') \otimes_{\mathbb{Q}[\pi_K]} \mathbb{Q}[\pi_K][X]/(X^n-\pi_K) \longrightarrow \text{\rm Hom}^{\! \circ \!}(W,W) \simeq \\[0.5 ex]
\text{\rm Hom}^{\! \circ \!}(\prod_{i=0}^{n-1} \sigma_{K|k}^{-i}(A'),A') \simeq \bigoplus_{i=0}^{n-1} \text{\rm Hom}^{\! \circ \!}(\sigma_{K|k}^{-i}(A'),A')
\end{array}
\end{equation}
is an isomorphism. Since we already know the homomorphism in the theorem to be a ring-homomorphism, this will conclude the proof.

Let ${\sigma_k} \in \text{Gal}(\overline{k}|k)$ be the Frobenius automorphism. By base-change, this induces an automorphism $\sigma_k$ of the $k$-scheme $W_{\overline{k}}$.

The endomorphism $\pi_k : W \longrightarrow W$ is uniquely determined by the fact that it operates on $\overline{k}$-valued points $P$ of $W_{\overline{k}}$ as the inverse of the ``arithmetic Frobenius operation'': $(\pi_k \otimes_k \text{id}_{\overline{k}}) \circ P = \sigma_k^{-1}(P)$.

Let $P = (P_i)_{i=0}^{n-1}$ be a $\overline{k}$-valued point of ${W}_{\overline{k}} \simeq \prod_{i=0}^{n-1} \sigma_{K|k}^{-i}(A')_{\overline{K}}$. Then $\sigma_k^{-1}(P) = (\sigma_k^{-1}(P_{i-1}))_{i=0}^{n-1}$ (where $P_{-1} := P_{n-1}$); see Subsection \ref{construction}. Thus $(\pi_k \otimes_k \text{id}_{\overline{k}}) \circ P = \sigma_k^{-1}(P) = (\sigma_k^{-1}(P_{i-1}))_{i=0}^{n-1} = ((\pi_k \otimes_k \text{id}_{\overline{k}}) \circ P_{i-1})_{i=0}^{n-1}$.

It follows that under the isomorphism $W_K \simeq \prod_{i=0}^{n-1} \sigma_{K|k}^{-i}(A')$, the endomorphism $\pi_k \otimes_k \text{id}_K$ of $W_K$ is given by the ``matrix''
\[ \left( \begin{array}{cccc}
0     & \cdots & \cdots & \pi_k \\
\pi_k & 0      & \cdots & 0 \\
0     & \ddots & \ddots & \vdots \\
0     & \ddots & \pi_k  & 0
\end{array} \right). \]
For $\lambda \in \text{\rm End}^{\! \circ \!}(A')$, $\Res^K_k(\lambda) \otimes_k \text{id}_K$ is given by the diagonal ``matrix''
\[ \left( \begin{array}{cccc}
\lambda     &  &  &  \\
& \sigma_{K|k}^{-1}(\lambda) & & \\
& & \ddots& \\
& & & \sigma_{K|k}^{-(n-1)}(\lambda)
\end{array} \right); \]
see Subsection \ref{functor-Res}. Let $x$ denote the image of $X$ in $\mathbb{Q}[\pi_K][X]/(X^n-\pi_K)$. Let $\lambda_1 x + \lambda_2 x^2 + \cdots + \lambda_{n} x^{n} \in \text{\rm Hom}^{\! \circ \!}(A',A') \otimes_{\mathbb{Q}[\pi_K]} \mathbb{Q}[\pi_K][X]/(X^n-\pi_K)$ where $\lambda_i \in \text{\rm End}^{\! \circ \!}(A')$. Such an element is mapped under the homomorphism of the theorem to an endomorphism of $W$ which is represented by the ``matrix''
{\small  \[ \left( \begin{array}{ccccc}
\lambda_n \, \pi_k^n & \lambda_{n-1} \, \pi_k^{n-1} & \cdots & \lambda_2 \, \pi_k^2 & \lambda_1 \, \pi_k \\[0.3 ex]
\sigma_{K|k}^{-1}(\lambda_{1}) \, \pi_k & \sigma_{K|k}^{-1}(\lambda_n) \, \pi_k^n & & \sigma_{K|k}^{-1}(\lambda_3) \, \pi_k^3 & \sigma_{K|k}^{-1}(\lambda_2) \, \pi_k^2 \\[0.3 ex]
\vdots & & \ddots & & \vdots \\[0.3 ex]
\sigma_{K|k}^{2-n}(\lambda_{n-2}) \, \pi_k^{n-2} & \sigma_{K|k}^{2-n}(\lambda_{n-3}) \, \pi_k^{n-3} & & \sigma_{K|k}^{2-n}(\lambda_{n}) \, \pi_k^n & \sigma_{K|k}^{2-n}(\lambda_{n-1}) \, \pi_k^{n-1} \\[0.3 ex]
\sigma_{K|k}^{1-n}(\lambda_{n-1}) \, \pi_k^{n-1} & \sigma_{K|k}^{1-n}(\lambda_{n-2}) \, \pi_k^{n-2}  & \cdots & \sigma_{K|k}^{1-n}(\lambda_1) \, \pi_k & \sigma_{K|k}^{1-n}(\lambda_n) \, \pi_k^n
\end{array} \right). \]}

The elements of $\text{\rm Hom}^{\! \circ \!}(A',A') \otimes_{\mathbb{Q}[\pi_K]} \mathbb{Q}[\pi_K][X]/(X^n-\pi_K)$ have a unique representation as $\lambda_1 x + \lambda_2 x^2 + \cdots + \lambda_{n} x^{n}$ with $\lambda_i \in \text{\rm End}^{\! \circ \!}(A')$. Under (\ref{finite-tensor-iso-row-vector}), this element corresponds to the first row in the above matrix, i.e.\ to the row vector 
\[( \begin{array}{cccc}
\lambda_n \pi_k^n & \lambda_{n-1} \pi^{n-1} & \cdots & \lambda_1 \pi_k \end{array}).\]
Now, every element of $\bigoplus_{i=0}^{n-1} \text{\rm Hom}^{\! \circ \!}(\sigma_{K|k}^{-i}(A'),A')$ has this form with unique $\lambda_i$. Thus (\ref{finite-tensor-iso-row-vector}) is an isomorphism.
\qed

\begin{remark}
Since the Frobenius endomorphism has degree a power of $p=\text{char}(k)$, we obtain in fact an isomorphism
\[ (\text{\rm End}(A')\otimes_{\mathbb{Z}[\pi_K]} \mathbb{Z}[\pi_K][X]/(X^n-\pi_K))\otimes\mathbb{Z}[1/p]\longrightarrow \text{\rm End}(W)\otimes\mathbb{Z}[1/p]. \]
\end{remark}

\begin{corollary}
\label{end-res-commutative}
$\text{\rm End}^{\! \circ \!}(W)$ is commutative if and only if $\text{\rm End}^{\! \circ \!}(A')$ is commutative.
\end{corollary}
The isomorphism of Theorem \ref{finite-end-structure} implies that the corresponding centers are isomorphic. Recalling from \cite{Ta1} that $Z(\text{\rm End}^{\!\circ\!}(A')) = \mathbb{Q}[\pi_K]$, we thus get:
\begin{corollary}
\label{center-finite-end-structure}
We have an isomorphism $\mathbb{Q}[\pi_K][X]/(X^n - \pi_K) \simeq Z(\text{\rm End}^{\! \circ \!}(W))$.
\end{corollary}

\subsection{Simplicity of the Weil restriction}
We are interested in the question whether the Weil restriction $W$ is simple.

In order that $W$ be simple, it is obviously necessary that $A'$ is simple. Furthermore, it is necessary that $A'$ is not isogenous to any abelian $K$-variety which can be defined over any proper intermediate field $\lambda$ of $K|k$ (i.e.\ any field $\lambda$ with $k \subseteq \lambda \subsetneq K$). (This holds for arbitrary finite separable field extensions $K|k$.)

For assume that this is the case. Since the scalar restriction of an isogeny is an isogeny, we can assume that $A'$ itself can be defined over such a $\lambda$; $A' = A_\lambda$ for some $\lambda$ as above and an abelian $\lambda$-variety $A$. By the defining functorial property of $W= \Res_k^K(A')$, we have a canonical homomorphism $\Res_k^\lambda(A) \longrightarrow W$ which is easily seen to be an immersion. Since the dimension of the immersed abelian variety is strictly smaller, $W$ is not simple.\\

We thus make the following assumption:

\emph{$A'$ is a simple abelian $K$-variety which is not isogenous to any abelian $K$-variety which can be defined over some proper intermediate field $\lambda$ of $K|k$.}

\begin{lemma}
\label{p_K-notin}
Under our assumption on $A'$, there does not exist a divisor $q$ of $n$ $(q \neq 1)$ such that $\pi_K \in \mathbb{Q}[\pi_K]^q$.
\end{lemma}
\emph{Proof} Assume that such a $q$ exists and let $\beta \in \mathbb{Q}[\pi_K]$ be such that $\beta^q = \pi_K$. (In particular $\mathbb{Q}[\pi_K] = \mathbb{Q}[\beta]$.)

Let $\lambda$ be the subfield of $K|k$ of index $q$, let $V$ be the Weil restriction of $A'$ with respect to $K|\lambda$. Denoting by $\chi$ characteristic polynomials of Frobenius-actions on Tate-modules we have $\chi_{V}(T) = \chi_{A'}(T^q)$, and $\beta$ is a root of $\chi_{V}$. This follows from the well-known fact that the operation of the absolute Galois group of $\lambda$ on $V(\overline{K})$ is induced by the operation of the absolute Galois group of $K$ on $A'(\overline{K})$; see \cite[\S 1,a)]{Mi-AA}.

It is easy to see that $V$ contains a simple abelian $\lambda$-variety $A$ such that the characteristic polynomial of the Frobenius of $A$ has $\beta$ as a root.

The structure of the endomorphism algebra $\text{End}^{\! \circ \!}(A)$ can be calculated from $\mathbb{Q}[\beta]$ as abstract field with generator $\beta$; see Subsection \ref{Honda-Tate-results}. Inserting $\beta$ and $\pi_K$ into formula (\ref{invariants}), one sees that the central-simple $\mathbb{Q}[\pi_K]$-algebras $\text{End}^{\! \circ \!}(A)$ and $\text{End}^{\! \circ \!}(A')$ have the same local invariants, thus they are isomorphic. Since by formula (\ref{dimension-from-invariants}), the dimension of abelian varieties can be calculated from their endomorphism algebras, it follows that $\text{dim}(A) = \text{dim}(A')$.

The immersion $A \longrightarrow V = \Res_\lambda^K(A')$ induces by the defining functorial property of the Weil restriction a non-trivial homomorphism $A_K \longrightarrow A'$. Since the dimensions agree and $A'$ is simple, this is an isogeny. A contradiction.
\qed

We now make use of the following well-known fact from field theory; see \cite[VI, \S 9, especially Theorem 9.1]{La}:
\begin{lemma}
Let $F$ be a field, $\alpha \in F, \alpha \neq 0$ and $n \in \mathbb{N}$. Assume that $\alpha \notin F^q$ for all prime divisors $q$ of $n$. Then
either $X^n-\alpha$ is irreducible over $F$ or $4 | n$ and $\alpha \in -4F^4$.
\end{lemma}

Together with Corollary \ref{center-finite-end-structure}, this implies:
\begin{proposition}
\label{finite-isotypic}
Under our assumption on $A'$, 
\begin{itemize}
\item
either $\Res^K_k(A')$ has exactly one isotypic component, i.e.\ all simple abelian subvarieties are isogenous
\item
or $4 | n$ and $\pi_K \in -4 \mathbb{Q}[\pi_K]^4$.
\end{itemize}
\end{proposition}
\emph{Proof}
By the previous two lemmata, under our assumption on $A'$,
either $X^n - \pi_K$ is irreducible over $\mathbb{Q}[\pi_K]$ 
or $4 |n$ and $\pi_K \in -4 \mathbb{Q}[\pi_K]^4$. Corollary \ref{center-finite-end-structure} implies: $X^n - \pi_K$ is irreducible over $\mathbb{Q}[\pi_K]$ if and only if $Z(\text{End}^{\! \circ \!}(\Res_k^K(A')))$ is a field. This in turn is equivalent to the fact that $\Res^K_k(A')$ has exactly one isotypic component.
\qed

\begin{remark}
By Honda's Theorem (see Proposition \ref{Honda}), it is obviously possible that
additionally to our general assumption on $A'$ the second condition is satisfied. It is interesting to note that there even exist ordinary elliptic curves $E'$ over fields of the form $\mathbb{F}_{p^4}$ ($p$ prime) which are non-isogenous to any elliptic $\mathbb{F}_{p^4}$-curve which can be defined over $\mathbb{F}_{p^2}$ and which satisfy $\pi_K \in -4\mathbb{Q}[\pi_K]^4$. Then $\Res^{\mathbb{F}_{p^4}}_{\mathbb{F}_p}(E')$ has more than one isotypic component. Since on the other hand it cannot contain an elliptic curve by our first assumption on $E'$, $\Res^{\mathbb{F}_{p^4}}_{\mathbb{F}_p}(E')$ has exactly two isotypic components both of which are simple.

For example, let $p$ be a prime such that $\left(\frac{-2}{p} \right) = 1$, let $K:=\mathbb{F}_{p^4}, k:= \mathbb{F}_p$.

By assumption, $p$ splits in the field $\mathbb{Q}[\sqrt{-2}]$; see \cite[Satz 8.5.]{Ne}. Since this field has class number 1, there is a prime element $\nu \in \mathcal{O}_{\mathbb{Q}[\sqrt{-2}]}$ such that $(\nu) (\overline{\nu}) = (p)$. (Where $-$ denotes conjugation.) Since the norm of an element is always positive, this implies $\nu \overline{\nu} = p$. If $i \in \mathbb{N}$, then $\nu^i \neq \overline{\nu}^i$, thus $\nu^i \notin \mathbb{Q}$.

Let $\alpha:=-\nu^4$. Then $\alpha^i \notin \mathbb{Q}$ for all $i \in \mathbb{N}$. In particular, $\mathbb{Q}[\alpha] = \mathbb{Q}[\sqrt{-2}]$. Let $E'$ be a simple abelian $K$-variety which corresponds to $(\mathbb{Q}[\alpha],\alpha)$ by Honda's Theorem (see Proposition \ref{Honda}). By formula (\ref{invariants}), all local invariants of $\text{End}^{\! \circ \!}(E')$ are congruent to 0, thus $\text{End}^{\! \circ \!}(E') \simeq \mathbb{Q}[\alpha]$, and $E'$ is an elliptic $K$-curve. Since $\alpha^i \notin \mathbb{Q}$ for all $i \in \mathbb{N}$, $E'$ is ordinary.

The algebraic integer $\alpha = -\nu^4 = -4 (\frac{\nu}{\sqrt{-2}})^4$ lies in $-4 \mathbb{Q}[\alpha]^4$. It remains to check that $E'$ is not isogenous to any elliptic $K$-curve which can be defined over $\mathbb{F}_{p^2}$.

Assume this was the case. Then there is a $\beta \in \text{End}^{\! \circ \!}(E') = \mathbb{Q}[\alpha] = \mathbb{Q}[\sqrt{-2}]$ with $\beta^2 =\alpha = - \nu^4$. This implies $i =\sqrt{-1} = \frac{\beta}{\nu^2} \in \mathbb{Q}[\sqrt{-2}]$, a contradiction.
\end{remark}
Our aim is now to give conditions under which the Weil restriction of $A'$ is even simple.

\begin{theorem}
\label{new-simple-res-simple}
Let $K|k$ be an extension of finite fields of degree $n$ and $A'$ a simple abelian $K$-variety. Assume that $A'$ is not isogenous to any abelian $K$-variety which can be defined over a proper intermediate field of $K|k$. Assume in addition that one of the following holds:
\begin{itemize}
\item
$\text{\rm End}(A')$ is commutative and further, if $4 | n$, then $\pi_K \notin -4 \mathbb{Q}[\pi_K]^4$.
\item
The extension degree $n$ is prime.
\end{itemize}
Then $\Res^K_k(A')$ is simple.
\end{theorem}

\emph{Proof}
Assume as in the theorem that $A'$ is not isogenous to any abelian variety which can be defined over a proper intermediate field of $K|k$.

We first treat the case that $\text{End}(A')$ is commutative and further, if $4 |n$, then $\pi_K \notin -4 \mathbb{Q}[\pi_K]^4$. Under these conditions, $\text{End}^{\! \circ \!}(\Res^K_k(A'))$ is also commutative (see Corollary \ref{end-res-commutative}), and by the above Proposition, $\Res^K_k(A')$ has exactly one isotypic component. This implies that $\Res^K_k(A')$ is simple.

We now come to the case that the extension degree $n$ is a prime. 
Let $B \subseteq \Res^K_k(A')$ be a simple abelian subvariety. Applying base-change, we get $B_K \subseteq \prod_{i=0}^{n-1} \sigma_{K|k}^{-i}(A')$. This implies dim$(A') \, | \, \text{dim}(B)$. Additionally, the dimensions cannot be equal since otherwise by the defining functorial property of the Weil restriction, we would have an isogeny $B_K \longrightarrow A'$ which is impossible by assumption. On the other hand, since by Proposition \ref{finite-isotypic} $\Res^K_k(A')$ has exactly one isotypic component, $\text{dim}(B) \, | \, \text{dim}(\Res^K_k(A')) = n \, \text{dim}(A')$. Since $n$ is a prime, this implies $\text{dim}(B) = \text{dim}(\Res^K_k(A'))$ thus $B = \Res^K_k(A')$.
\qed

\begin{remark}
Let $K:=\mathbb{F}_{p^4}, k:=\mathbb{F}_{p}$ where $p$ is a prime with $p \equiv 1 \; (\text{mod } 4)$. We will now give an elliptic $K$-curve $E'$ with \emph{non-commutative} endomorphism ring such that $\Res^K_k(E')$ is \emph{non-simple} even though $E'$ is not isogenous to any abelian $\mathbb{F}_{p^4}$-variety which can be defined over $\mathbb{F}_{p^2}$ and the condition $\pi_K \notin -4 \mathbb{Q}[\pi_K]^4$ is satisfied.

Let $E'$ be a simple abelian $K$-variety which corresponds to the integer $-p^2$ by Honda's Theorem; see Proposition \ref{Honda}. By formula (\ref{invariants}), the local invariants of $\text{End}^{\! \circ \!}(E')$ at $p$ and $\infty$ are congruent to $\frac{1}{2}$, thus $E'$ is a super-singular elliptic curve such that all endomorphisms of $E'_{\overline{\mathbb{F}}_{\!p}}$ can be defined over $\mathbb{F}_{p^4}$.

Assume there is an elliptic $\lambda := \mathbb{F}_{p^2}$-curve $E$ such that $E_K \sim E'$. Let $\pi_\lambda$ be its Frobenius endomorphism. Then we have $\mathbb{Q}[\pi_\lambda] \simeq \mathbb{Q}[i]$ ($i :=\sqrt{-1}$), and under this isomorphism, $\pi_\lambda$ corresponds to $ip$. Now by assumption, $p$ splits in $\mathbb{Q}[i]$, and from formula (\ref{invariants}), it follows that the local invariants of $\text{End}^{\! \circ \!}(E)$ over $p$ are congruent to $\frac{1}{2}$, thus by (\ref{dimension-from-invariants}), $E$ is 2-dimensional, a contradiction.

Let $W$ be the Weil restriction of $E'$ with respect to $K|k$. Then by Corollary \ref{center-finite-end-structure}, the center of $\text{End}^{\! \circ \!}(W)$ is isomorphic to $\mathbb{Q}[X]/(X^4+p^2) = \mathbb{Q}[\sqrt[4]{-p^2}]$, and under this isomorphism $\pi_k$ corresponds to $\sqrt[4]{-p^2}$. In this field, $p$ is ramified of degree 2 and splits into 2 prime ideals (because it already splits in the subfield $\mathbb{Q}[i]$). Again by formula (\ref{invariants}), the endomorphism algebras of the simple components of $W$ are fields, thus isomorphic to $\mathbb{Q}[\sqrt[4]{-p^2}]$. It follows with (\ref{dimension-from-invariants}) that the simple components of $W$ are 2-dimensional, thus $W$ is not simple.
\end{remark}

\subsection{Appendix to Section 2: Some results by Honda and Tate}
\label{Honda-Tate-results}
For the convenience of the reader, we recall Honda's Theorem on the classification of simple abelian varieties over finite fields and Tate's results how to compute the structure of the endomorphism ring of an abelian variety over a finite field; c.f. \cite{Ho, Ta1, Ta2}.

Fix a finite field $k= \mathbb{F}_q$, where $q = p^a$ with $p$ a prime and $a \in \mathbb{N}$. Then, if $A$ is a simple abelian $k$-variety and $\pi_k$ is its Frobenius endomorphism, for every inclusion $\varphi$ of $\mathbb{Q}[\pi_k]$ into $\overline{\mathbb{Q}}$, we have $|\varphi(\pi_k)| = q^{\frac{1}{2}}$.

Now Honda's Theorem states:

\begin{proposition}[Honda]
\label{Honda}
The assignment $A \mapsto (\mathbb{Q}[\pi_k], \pi_k)$ induces a bijection between the set of isogeny classes of simple abelian $k$-varieties and the set of isomorphism classes of fields $\mathbb{Q}[\alpha]$ with fixed generator $\alpha$ such that $\alpha$ is an algebraic integer and under all inclusions into $\overline{\mathbb{Q}}$, $\alpha$ has absolute value $q^{\frac{1}{2}}$.
\end{proposition}

By Honda's Theorem, for every simple abelian $k$-variety $A$, the structure of the endomorphism algebra $\text{End}^{\! \circ \!}(A)$ only depends on $\mathbb{Q}[\pi_k]$ as abstract field with generator $\pi_k$. Since $\text{End}^0(A)$ is central-simple over $\mathbb{Q}[\pi_k]$, to determine its structure, we only have to give its local invariants at all finite and real valuations.

The formula for this is as follows:
Let $v$ be a normalized valuation of $\mathbb{Q}[\pi_k]$. Then, if $v$ is finite, the local invariant of $\text{End}^{\! \circ \!}(A)$ at $v$ is given by
\begin{equation}
\label{invariants}
\text{inv}_v \equiv \frac{v(\pi_k)}{a}f_v \; \; (\text{mod} \; 1),
\end{equation}
where $f_v$ denotes the absolute residue degree of $\mathbb{Q}[\pi_k]$ at $v$. In particular, if $v$ is a finite valuation which does not lie over the valuation of $p$, the local invariant is congruent to 0.

If $v$ is real, then the local invariant is congruent to $\frac{1}{2}$. 

Let $m$ be the least common denominator of the local invariants. Then the order of $\text{End}^{\! \circ \!}(A)$ in the Brauer group of $\mathbb{Q}[\pi_k]$ is $m$, $m^2 = [\text{End}^{\! \circ \!}(A) : \mathbb{Q}[\pi_k]]$, and the dimension of $A$ in given by
\begin{equation}
\label{dimension-from-invariants}
 \text{dim}(A) = \frac{1}{2} \, m \, [\mathbb{Q}[\pi_k]:\mathbb{Q}].
\end{equation}

\section{Results for abelian varieties which can be defined over the base-field}
Throughout this section, let $K|k$ be a finite Galois extension of degree $n$ with Galois group $G$, and let $A$ be an abelian $k$-variety of dimension $d$. Let $W$ be the Weil restriction of $A_K$ with respect to $K|k$.

We want to determine the structure of the endomorphism ring of $W$, and the isogeny decomposition of $W$ over $k$.

\subsection{Arithmetic becomes geometric operation}
\label{subsection-arith-geo}

For any $k$-scheme $Z$, $G$ operates on $A_K(Z_K)$ 
by $ \tau(P) = \tau P \tau^{-1}$. These operations define an automorphism of the functor $Z \mapsto A_K(Z_K)$ from the category of $k$-schemes to the category of abelian groups. 
We obtain automorphisms of the representing object $W=\Res^{K}_k(A_K)$ which we denote by $a_\tau$ for $\tau\in G^{\text{opp}}$. We thus have a group-homomorphism $a : G^{\text{opp}} \longrightarrow \text{Aut}(W), \, \tau \longrightarrow a_\tau$, where $\text{Aut}(W)$ denotes the group of automorphisms of the abelian $k$-variety $W$.

We want to calculate how $a_\tau \otimes_k \text{id}_K$ operates on $W_K \simeq A^{G^{\text{opp}}}_{K}$.

We have $\tau(u) = \tau(p_{\text{id}}) = p_{\tau^{-1}} : \, W_K \longrightarrow A_K$ by (\ref{operation-on-projection}). 
The homomorphism $a_\tau$ of the abelian $k$-variety $W$ is the $W$-valued point of $W$ which corresponds to $\tau(u)$. So by Subsection \ref{construction}, $a_\tau \otimes_k \text{id}_K = (\sigma^{-1}(\tau(u)))_{\sigma \in G^{\text{opp}}} = (\sigma^{-1}(p_{\tau^{-1}}))_{\sigma \in {G^{\text{opp}}}} = (p_{\tau^{-1} \sigma})_{\sigma \in {G^{\text{opp}}}}$. (The last equation follows from (\ref{operation-on-projection}).) We have established:

\begin{lemma}
\label{lemma-arith-geo}
$a_\tau \otimes_k \text{\rm id}_K : A^{G^{\text{\rm opp}}}_{K} \longrightarrow A^{G^{\text{\rm opp}}}_{K}$ operates on $Z$-valued points (any $Z$) by $(P_\sigma)_{\sigma \in {G^{\text{\rm opp}}}} \mapsto (P_{\tau^{-1} \sigma})_{\sigma \in {G^{\text{\rm opp}}}}$.
\end{lemma}

\label{factors}
\subsection{The endomorphism ring as skew group ring}

\begin{lemma}
\label{geo-arith-skew}
Let $\tau \in G^{\text{\rm opp}}, \lambda \in \text{\rm End}(A_K)$. Then $a_\tau \circ \Res^K_k(\lambda) = \Res^K_k(\tau(\lambda)) \circ a_\tau \in \text{\rm End}(W)$.
\end{lemma}
\emph{Proof} Easy calculation on $Z$-valued points.\qed

To formulate the result about the structure of the endomorphism ring of $W$, we need a generalization of the concept of a group ring first.

\paragraph{Definition}
Let $\Lambda$ be a ring, $G$ a group, $t : G \longrightarrow \text{\rm Aut}(\Lambda)$ a group-homo\-morphism. The application of $t(\sigma)$ to some $\lambda \in \Lambda$ will by denoted by $\sigma(\lambda)$. 
Following \cite{Pa}, we define the \emph{skew group ring} $\Lambda^t[G]$ to be the following ring:\footnote{This ring is a special case of a \emph{crossed product} (with respect to some operation); cf. \cite{Pa}. In \cite{CR}, the same ring is called \emph{twisted group ring}. However, in \cite{Pa}, this word is reserved for the special case of a crossed product with respect to a trivial group operation.}
The underlying abelian group is $\Lambda^G$ with the usual ``componentwise'' addition. 
As usual, for $\tau \in G$, let $\tau$ also denote $(\delta_{\sigma,\tau})_{\sigma \in G} \in \Lambda^G$.
The multiplication is defined by $\sum_{\sigma \in G} \lambda_\sigma \, \sigma \cdot \sum_{\nu \in G} \mu_\nu \, \nu = \sum_{\sigma, \nu \in G} \lambda_\sigma \, \sigma(\mu_\nu) \, \sigma \nu$.\\

The ring $\Lambda$ is naturally immersed in $\Lambda^t[G]$. For fixed $\Lambda, \, G$ and $t : G \longrightarrow \text{Aut}(\Lambda)$, the ring $\Lambda^t[G]$ has the following universal property:
\begin{lemma}
\label{skew-group-ring-universal-property}
Let $B$ be a ring, $f: \Lambda \longrightarrow B$ be a ring-homomorphism, and let $g : G \longrightarrow B^*$ be a group-homomorphism. Assume that for $\lambda \in \Lambda, \tau \in G$, $g(\tau) \, f(\lambda) = f(\tau(\lambda)) \, g(\tau)$. Then there is a unique ring-homomorphism $\Lambda^t[G] \longrightarrow B$ with $\lambda \mapsto f(\lambda)$ and $\tau \mapsto g(\tau)$.
\end{lemma}

Now let $G$ be the Galois group as above, $t: G^{\text{opp}} \longrightarrow \text{\rm Aut}(\text{\rm End}(A_K))$ the natural operation given by $\sigma \mapsto (\lambda \mapsto \sigma(\lambda) = \sigma \lambda \sigma^{-1})$. From Lemmata \ref{geo-arith-skew} and \ref{skew-group-ring-universal-property} it follows that $\sum_{\sigma \in G^{\text{opp}}} \lambda_\sigma \, \sigma \mapsto \sum_{\sigma \in G^{\text{opp}}} \Res^K_k(\lambda_\sigma) \, a_\sigma$ defines a ring-homomorphism
\begin{equation}
\label{skew-group-ring-end-W-iso}
\text{\rm End}(A_K)^t[G^{\text{\rm opp}}] \longrightarrow \text{\rm End}(W).
\end{equation}
\begin{theorem}
\label{theorem-skew-group-ring-end-W-iso}
Let $K|k$ be a finite Galois extension with Galois group $G$, $A$ an abelian $k$-variety, $W$ the Weil restriction of $A_K$ with respect to $K|k$, $t : G^{\text{\rm opp}} \longrightarrow \text{\rm Aut}(\text{\rm End}(A_K))$ the natural operation. Then
\[\text{\rm End}(A_K)^t[G^{\text{\rm opp}}] \longrightarrow \text{\rm End}(W),\, \sum_{\sigma \in G^{\text{\rm opp}}} \lambda_\sigma \, \sigma \mapsto \sum_{\sigma \in G^{\text{\rm opp}}} \Res^K_k(\lambda_\sigma) \, a_\sigma\]
 is an isomorphism.
\end{theorem}
\emph{Proof}
Analogously to the proof of Theorem \ref{finite-end-structure}, we make use of the isomorphism $\text{\rm Hom}(W,W) \simeq \text{\rm Hom}(A_K^{\; G^{\rm \text{opp}}},A_K) \simeq \bigoplus_{\sigma \in G^{\text{opp}}} \text{\rm Hom}(A_K,A_K)$ of the right-hand side.

By (\ref{operation-on-projection}), the image of some $\sigma \in G^{\text{opp}}$ in $\text{\rm Hom}(A^{\; G^{\text{opp}}}_K,A_K)$ is $p_{\sigma^{-1}}$, corresponding to the row vector which is zero except at the ``$\sigma$-th'' entry where it is \nolinebreak 1.

Thus the image of $\sum_{\sigma \in G^{\text{opp}}} \lambda_\sigma \sigma$ (where $\lambda_\sigma \in \text{\rm End}(A_K)$) is $\sum_{\sigma \in G^{\text{opp}}} \lambda_{\sigma^{-1}} \, p_\sigma$, corresponding to the row vector $(\lambda_{\sigma^{-1}})_{\sigma \in G^{\text{opp}}}$.

It is thus immediate that we have an isomorphism.
\qed

\begin{corollary}
\label{corollary-skew-group-ring-end-W-iso}
The isomorphism in the theorem induces an isomorphism \linebreak
$\text{\rm End}^{\! \circ \!}(A_K)^t[G^{\text{\rm opp}}] \longrightarrow \text{\rm End}^{\! \circ \!}(W)$.
\end{corollary}

By the Complete Reducibility Theorem (see \cite[Proposition 12.1]{Mi-2}) we know that the ring $\text{\rm End}^{\! \circ \!}(W)$ is semi-simple. Thus the skew group ring \linebreak ${\text{\rm End}^{\! \circ \!}(A_K)}^t[G{^\text{opp}}]$ is semi-simple.

It can be proven more generally that every crossed product over a semi-simple ring with a finite group in which the group order is invertible is semi-simple; see \cite[Theorem 4.1.]{Pa}.\\

We now want to study the ring-homomorphism
\begin{equation}
\label{end-right-regular-representation}
\begin{array}{c}
\text{\rm End}(A_K)^t[G^{\text{opp}}] \stackrel{\sim}{\longrightarrow} \text{\rm End}(W) \hookrightarrow \\[0.5 ex]
\text{\rm End}(W_K) \simeq \text{\rm End}(A_K^{\; G^{\text{opp}}}) \simeq \text{\rm M}_{G^{\text{opp}}}(\text{\rm End}(A_K)).
\end{array}
\end{equation}
We denote the matrix corresponding to $a_\tau$ by $A_\tau$ and the matrix corresponding to $\Res^K_k(\lambda)$ by $J(\lambda)$ (for
$a_{\tau}$ as above and $\lambda\in\text{End}(A_K)$).

We have already shown in Subsection \ref{functor-Res} that $J(\lambda)$ is the diagonal matrix $(\sigma^{-1}(\lambda) \delta_{\sigma,\nu})_{\sigma,\nu \in G^{\text{opp}}}$.

Let us determine to which matrix $A_{\tau}\in \text{\rm M}_{G^{\text{opp}}}(\text{End}(A_K))$ the endomorphism $a_\tau$ corresponds. First of all, $p_\sigma : W_K \simeq A_K^{\; G^{\text{opp}}} \longrightarrow A_K$ corresponds to the row vector $(\delta_{\sigma, \nu})_{\nu \in G^{\text{opp}}}$. As $a_\tau = (p_{\tau^{-1} \sigma})_{\sigma \in G^{\text{opp}}}$ (see Lemma \ref{lemma-arith-geo}), we get
\begin{equation}
\label{a-A-matrix}
A_\tau = (\delta_{\tau^{-1} \sigma, \nu})_{\sigma, \nu \in G^{\text{opp}}} = (\delta_{\sigma, \tau \nu})_{\sigma, \nu \in G^{\text{opp}}}.
\end{equation}

Before continuing let us recall the definition of the left regular (matrix) representation.

\subsubsection*{The left regular (matrix) representation}
Let $\Lambda$ be a ring. If $\Lambda \longrightarrow \Xi$ is a homomorphism of rings, we can regard $\Xi$ as $\Lambda$-right module, and if we do so, we write $\text{\rm End}^r_\Lambda(\Xi)$ for the ring of endomorphisms.

Now let $\Lambda \longrightarrow \Xi$ be a homomorphism of rings and assume additionally that $\Xi$ is free as $\Lambda$-right module on a finite set of generators $\Sigma$, i.e.\ $\Xi \simeq \Lambda^\Sigma$ as $\Lambda$-right modules. Multiplication by elements of $\Xi$ from the left induces a ring-homomorphism
\begin{equation}
\label{left-representation-abstract}
l : \Xi \longrightarrow \text{\rm End}^r_\Lambda(\Xi) \simeq \text{\rm End}^r_\Lambda(\Lambda^\Sigma),
\end{equation}
the \emph{left regular representation}. 

For a fixed basis $\Sigma$, the right-hand side of (\ref{left-representation-abstract}) is canonically isomorphic to the matrix ring $\text{\rm M}_\Sigma(\Lambda)$. The isomorphism is given as follows:
\begin{equation}
\label{Lambda-Sigma-end-matrix-ring}
\begin{array}{c}
\text{\rm End}^r_\Lambda(\Lambda^\Sigma) \longrightarrow \text{\rm M}_\Sigma(\Lambda), \; a \mapsto (\alpha_{\sigma,\nu})_{\sigma,\nu \in \Sigma} \text{ with } \alpha_{\sigma,\nu} \in \Lambda \\
\text{and } a(\nu) = \sum_{\sigma \in \Sigma} \sigma \, \alpha_{\sigma,\nu}.
\end{array}
\end{equation}
By composition of (\ref{left-representation-abstract}) with (\ref{Lambda-Sigma-end-matrix-ring}), we get the \emph{left regular matrix representation} (with respect to the basis $\Sigma$).
\[L : \Xi \longrightarrow \text{\rm M}_\Sigma(\Lambda).\]

We now apply these concepts in the context of the skew group ring. Let $G$ be a finite group, $t : G \longrightarrow \text{\rm Aut}(\Lambda)$ be a homomorphism, $\Lambda^t[G]$ the corresponding skew group ring.

We calculate explicitly the left regular representation $l :  \Lambda^t[G] \longrightarrow $\linebreak$\text{\rm End}_\Lambda^r(\Lambda^t[G])$ and the left regular matrix representation $L : \Lambda^t[G] \longrightarrow \text{\rm M}_G(\Lambda)$ with respect to the basis $G$.

Let $\tau \in G$. Then $l(\tau) : \nu \mapsto \tau \nu = \sum_{\sigma \in G} \sigma \delta_{\sigma,\tau \nu}$ and thus
\[L(\tau) = (\delta_{\sigma, \tau \nu})_{\sigma, \nu \in G}. \]

Let $\lambda \in \Lambda$. Then $l(\lambda) : \nu \mapsto \lambda \, \nu = \nu \, \nu^{-1}(\lambda)$ and thus
\[L(\lambda) = (\sigma^{-1}(\lambda) \, \delta_{\sigma,\nu})_{\sigma,\nu \in G}. \]
\bigskip

We are now going to relate these definitions and calculations with our situation. So let $\Lambda := \text{\rm End}(A_K)$, $G$ the Galois group and $t :G^{\text{\rm opp}} \longrightarrow \text{\rm End}(A_K)$ the natural operation. Let $L$ be the left regular matrix representation of $\Lambda$ with respect to the basis $G^{\text{opp}}$. Then $L(\tau) = A_\tau$ and $L(\lambda) = J(\lambda)$.
Thus:
\begin{proposition}
\label{Hom-is-left-regular-m.-r.}
Homomorphism (\ref{end-right-regular-representation}) is the left regular matrix representation of the skew group ring $\text{\rm End}(A_K)^t[G^{\text{\rm opp}}]$ with respect to the basis $G^{\text{\rm opp}}$.
\end{proposition}

\subsection{The Rosati involution}
Let $\varphi : A_K \longrightarrow \widehat{A}_K$ be a polarization. Then $\Res^K_k(\varphi) : W \longrightarrow \widehat{W}$ is also a polarization; see Subsection \ref{Weil restriction of pol. abelian varieties}.

We want to calculate how the Rosati involution of $W$ with respect to $\Res^K_k(\varphi)$ is given under the isomorphism of Corollary \ref{corollary-skew-group-ring-end-W-iso}.

Let us denote the Rosati involution by $(\ldots)'$.

First of all, the (defining) equation $\lambda' = \varphi^{-1} \widehat{\lambda} \varphi$ where $\lambda \in \text{\rm End}^{\! \circ \!}(A_K)$ implies
\[ \Res^K_k(\lambda') = \Res^K_k(\varphi)^{-1} \circ \Res^K_k(\widehat{\lambda}) \circ \Res^K_k(\varphi) = \Res^K_k(\lambda)'. \]
(This holds more generally for any abelian $K$-variety $A'$ instead of $A_K$.)

We use the inclusion of $\text{\rm End}^{\! \circ \!}(W)$ into the matrix ring $\text{\rm M}_{G^{\text{opp}}}(\text{\rm End}^{\! \circ \!}(A))$ and the fact that $\Res^K_k(\varphi) \otimes_k \text{id}_K$ is a product polarization to calculate the Rosati involution of $a_{\tau}$ with the help of Lemma \ref{matrix-Rosati}.

Since $a_\tau$ corresponds to the matrix $A_\tau = (\delta_{\sigma, \tau \nu})_{\sigma, \nu \in G^{\text{opp}}}$ (see (\ref{a-A-matrix})), $a_\tau'$ corresponds to the matrix $(\delta_{\nu, \tau \sigma})_{\sigma, \nu \in G^{\text{opp}}} = (\delta_{\tau^{-1} \nu, \sigma})_{\sigma, \nu \in G^{\text{opp}}} = $\linebreak$(\delta_{\sigma, \tau^{-1} \nu})_{\sigma, \nu \in G^{\text{opp}}} = A_{\tau^{-1}}$. Thus \[a_\tau' = a_{\tau^{-1}}.\]
Since the Rosati involution is an anti-ring-endomorphism, this implies:
\begin{proposition}
\label{old-Res-Rosati}
Let $K|k$ be a finite Galois field extension with Galois group $G$, $A$ an abelian $k$-variety, $W$ the Weil restriction of $A_K$ with respect to $K|k$. Let $\varphi: A \longrightarrow \widehat{A}$ be a polarization. Let $\lambda \mapsto \lambda'$ be the Rosati involution associated to $\varphi$. Then under the isomorphism of Corollary \ref{corollary-skew-group-ring-end-W-iso}, the Rosati involution associated to the polarization $\Res^K_k(\varphi) : W \longrightarrow \widehat{W}$ is given by $\sum_{\sigma \in G^{\text{\rm opp}}} \lambda_\sigma \, \sigma \mapsto \sum_{\sigma \in G^{\text{\rm opp}}} \sigma^{-1} \lambda'_\sigma = \sum_{\sigma \in G^{\text{\rm opp}}}  \sigma^{-1}(\lambda'_\sigma) \, \sigma^{-1}$.
\end{proposition}

\subsection{Dimensions of components}

As in the above proposition, let $A$ be an abelian $k$-variety, $K|k$ a galois field extension of degree $n$ with galois group $G$, $W$ the Weil restriction of $A_K$ with respect to $K|k$, and let $t : G^{\text{opp}} \longrightarrow \text{End}(A_K)$ be the natural operation.

Assume $D \subseteq \text{\rm End}^{\! \circ \!}(A_K)$ is a skew field, invariant under the operation $t$.

Let $\bigoplus_{i=1}^s \Lambda_i = D^t[G^{\text{opp}}]$ be a decomposition of the $D^t[G^{\text{opp}}]$-right module $D^t[G^{\text{opp}}]$. This defines a decomposition $1 = \sum_i e_i$ where the $e_i$ are orthogonal idempotents, $e_i \in \Lambda_i$, such that $\Lambda_i = e_i \, D^t[G^{\text{opp}}]$. Conversely, if we are given a decomposition $1 =\sum_i e_i$ with orthogonal idempotents $e_i$, then the $\Lambda_i := e_i \, D^t[G^{\text{opp}}]$ define a direct sum decomposition of the $D^t[G^{\text{opp}}]$-right module $D^t[G^{\text{opp}}]$.

Via the inclusion $D^t[G^{\text{opp}}] \hookrightarrow \text{\rm End}^{\! \circ \!}(A_K)^t[G^{\text{opp}}] \simeq \text{\rm End}^{\! \circ \!}(W)$, we can regard the $e_i$ to be elements of $\text{\rm End}^{\! \circ \!}(W)$. For each $i$, let $c_i \in \mathbb{N}$ such that $c_i e_i \in \text{\rm End}(W)$.

Now put $W_i := (c_i e_i)(W)$. The $W_i$ are abelian subvarieties of $W$ and $\bigoplus_{i=1}^s W_i \sim W$. (Conversely, such an isogeny decomposition where the $W_i$ are abelian subvarieties of $W$ determines a decomposition of $\text{\rm End}^{\! \circ \!}(W)$ as right-$\text{\rm End}^{\! \circ \!}(W)$ module.)

\begin{proposition}
\label{isotyic-components}
Let $D \subseteq \text{\rm End}^{\! \circ \!}(A_K)$ be a skew field, invariant under the operation $t$ on $\text{\rm End}^{\! \circ \!}(A_K)$. Let $\bigoplus_{i=1}^s \Lambda_i = D^t[G^{\text{\rm opp}}]$ be a decomposition of the $D^t[G^{\text{\rm opp}}]$-right module $D^t[G^{\text{\rm opp}}]$. This corresponds to a decomposition $\text{\rm id}_{A_K} = \sum_i e_i$. Let $W_i := (c_i e_i)(W)$ be as above. Then ${W_i}_K \approx A_K^{n_i}$ (non-canonical isomorphism) where
\[ n_i = \text{\rm dim}_D(\Lambda_i). \]
\end{proposition}
\emph{Proof}
Choose a bijection between $G^{\text{opp}}$ and the set $\{1, \ldots, n \}$. Then $A_K^{\; G^{\text{opp}}} \simeq A_K^n$.

Let $l$ and $L$ be the left regular (matrix) representations of $\text{\rm End}^{\! \circ \!}(A_K)^t[G^{\text{opp}}]$, $l_D$ and $L_D$ the left regular (matrix) representations of $D^t[G^{\text{opp}}]$ (both regular matrix representations with respect to the basis $G^{\text{opp}}$). Let $\iota_{\text{\rm M}} : \text{\rm M}_{G^{\text{opp}}}(D) \longrightarrow \text{\rm M}_{G^{\text{opp}}}(\text{\rm End}(A_K))$ be the canonical inclusion. Then $L = \iota_{\text{\rm M}} \, L_D$.

By construction $l_D(e_i)$ is the identity on $\Lambda_i$ and zero on all $\Lambda_j$ for $j \neq i$.

Let $n_i$ be the dimension of the $D$-module $\Lambda_i$. For each $i$, choose a basis $(b_i^{(j)})_{j=1, \ldots, n_i}$ of the $D$-module $\Lambda_i$. Then all $n$ elements $b_i^{(j)}$ define a basis of the $D$-module $D^t[G^{\text{opp}}]$. With respect to this basis, the matrix associated to $l_D(e_i)$ is zero outside a block of size $n_i$ where it is the identity matrix.

We now have two matrix representations of $l_D(e_i)$ with respect to different bases, and via a base change matrix, we can transform one into the other: There exists an invertible matrix $B\in \text{Gl}_n(D)$ such that $B L_D(e_i) B^{-1}$ is zero outside a block of size $n_i$ where it is the identity matrix. 

Let $b \in \text{End}^{\!\circ\!}(A_K^{G^{\text{opp}}}) \simeq \text{End}^{\!\circ\!}(A_K^n)$ correspond to $\iota_{\text{\rm M}}(B)$. By Proposition \ref{Hom-is-left-regular-m.-r.} and our notational conventions, the endomorphism associated to the matrix $L(e_i) = \iota_{\text{\rm M}} L_D(e_i)$ is $e_i \otimes_k \text{id}_K$. By the above considerations, $b(e_i \otimes_k \text{id}_K) b^{-1}$ is an endomorphism whose image is isomorphic to $A_K^{n_i}$. It follows that the image of $c_i \, e_i \otimes_k \text{id}_K$ is also isomorphic to $A_K^{n_i}$.
\qed

\begin{remark}
Let $A_K$ be simple, $D= \text{\rm End}^{\! \circ \!}(A_K)$. Assume that all $e_i$ in the above proposition are central. Then all $\Lambda_i$ as above are rings and we have an isomorphism $\prod_{i=1}^s \Lambda_i \simeq D^t[G^{\text{opp}}]\simeq \text{End}^{\! \circ \!}(W)$ of rings. Furthermore, the $(c_i e_i)(W)$ are generated by isotypic components of $W$ and $\Lambda_i \simeq \text{\rm End}^{\! \circ \!}((c_i e_i)(W))$. So in particular, the number $n_i$ in the above proposition satisfies $n_i = \text{dim}_D(\text{\rm End}^{\! \circ \!}(W_i))$.
\end{remark}

\subsection{The cyclic case}
\label{decomp-cyclic}
We now apply the above results to the case that $G$ is cyclic of order $n$.

We identify $G$ with $G^{\text{opp}}$ and fix some generator $\sigma \in G$. Let $a = a_\sigma \in \text{\rm End}(W)$ be the automorphism corresponding to $\sigma$.

Denote the residue class of $X$ in $\mathbb{Q}[X]/(X^n-1)$ by $x$. Then we have an inclusion
\[ \mathbb{Q}[X]/(X^n-1)  \longrightarrow  \text{\rm End}^{\! \circ \!}(A_K)^t[G], \; x \mapsto \sigma. \]
The polynomial $X^n -1 \in \mathbb{Z}[X]$ splits over $\mathbb{Z}$ as
\[ X^n -1 = \prod_{d|n} \Phi_d, \]
where $\Phi_d$ is the $d$-th cyclotomic polynomial.

Let $\Phi'_d := (X^n-1)/\Phi_d$. By the Euclidian algorithm, there exist $\Psi_d \in \mathbb{Q}[X]$ with $\sum_{d|n} \Psi_d \Phi'_d =1$. Let $E_d := \Psi_d \, \Phi'_d$. Then the $E_d(x) \in \mathbb{Q}[X]/(X^n-1)$ are pair-wise orthogonal idempotents. The corresponding decomposition is

\[ \mathbb{Q}[X]/(X^n-1) \simeq  \prod_{d|n} \mathbb{Q}[X]/\Phi_d = \prod_{d|n} \mathbb{Q}(\zeta_d). \]

(This is nothing but the Chinese Remainder Theorem in this particular case.)

Let $W_d := c_d \, E_d(a)(W)$ for suitable $c_d \in \mathbb{N}$. We then have an isogeny decomposition
\[ W \sim \prod_{d|n} W_d, \]
and by Proposition \ref{isotyic-components}, the $W_d$ are abelian varieties with ${W_d}_K \approx A_K^{\varphi(d)}$.

We also have $W_d = \Phi'_d(a)(W)$. -- We only have to show that $c_d \Phi'_d(a)(W) \subseteq W_d$. This follows from $\Phi'_d(x) = (\sum_{f|n} \Psi_f(x) \Phi'_f(x)) \Phi'_d(x) = \Psi_d(x) {\Phi'_d}^2(x) = E_d(x) \Phi'_d(x)$.

It is clear that $W_d$ is also the reduced identity component of the kernel of \[ \begin{array}{c}
c_d(\text{id} - E_d(a)) = c_d \sum_{f|n, f \neq d} \Psi_f(a) \Phi'_f(a) = \\
(c_d \sum_{f|n, \, f \neq d} \Psi_f(a) \prod_{g|n, \, g \neq d,f} \Phi_g(a)) \, \Phi_d(a).
\end{array} \]
It is also the reduced identity component of the kernel of $\Phi_d(a)$. -- We only have to show that $W_d$ is contained in this kernel. But since $W_d = \Phi'_d(a)(W)$ and $\Phi'_d(x) \Phi_d(x) = 0$, this is obvious.\\

We now want to study whether the $W_d$ are simple or split further. We make the following assumptions.

\emph{$A_K$ is a simple abelian $K$-variety whose endomorphisms can be defined over $k$ and whose endomorphism ring is commutative.}

Note that if $k$ is finite, all endomorphisms of $A_K$ can automatically be defined over $k$ if we assume $\text{\rm End}(A_K)$ to be commutative.

Also if $A$ is an ordinary elliptic curve over an arbitrary field $k$ and $n$ is odd, then all endomorphisms of $A_K$ can be defined over $k$. This is because under this condition, End($A_K$) is either $\mathbb{Z}$ or a quadratic order, thus the only possible non-trivial automorphism of $\text{\rm End}(A_K)$ has order 2, and consequently the image of the representation $\text{Gal}(K|k) \longrightarrow \text{\rm Aut}(\text{\rm End}(A_K))$ is trivial.

Under the assumptions, we have the isomorphisms
\[ \begin{array}{ccccc}
\text{\rm End}^{\! \circ \!}(A)[X]/(X^n-1) & \simeq &  \text{\rm End}^{\! \circ \!}(A_K)[G] & \simeq & \text{\rm End}^{\! \circ \!}(W)  \\
x & \mapsto & \sigma & \mapsto & a
\end{array}. \]
Let $\Phi_d$ split into the product of the non-trivial monic irreducible polynomials $\Phi_d^{(1)}, \Phi_d^{(2)}, \ldots, \Phi_d^{(r_{\!d})}$ over $\text{\rm End}^{\! \circ \!}(A)$. Let ${\Phi_{d}'}^{(i)} := (X^n-1)/\Phi_d^{(i)}$. Since $X^n-1$ is separable in characteristic zero, the $\Phi_d^{(i)}$ are pair-wise coprime for varying $d$ and $i$, and there exist $\Psi_d^{(i)}$ with $\sum_{d|n} \sum_{1=i}^{r_{\!d}} \Psi_d^{(i)} {\Phi'_d}^{(i)} =1$. Let $E_d^{(i)} := \Psi_d^{(i)} {\Phi'_d}^{(i)}$.

Let $W_d^{(i)} := c_d^{(i)} E_d^{(i)}(a)(W)$ for suitable $c_d^{(i)} \in \mathbb{N}$. Then again by Proposition \ref{isotyic-components}, $W_d^{(i)}$ is an abelian $k$-variety with $(W_{d}^{(i)})_K \approx A_K^{\text{deg}(\Phi_d^{(i)})}$. The abelian $k$-variety $W_{d}^{(i)}$ is simple and its endomorphism algebra is isomorphic to the field $\text{\rm End}^{\! \circ \!}(A)[X]/\Phi_d^{(i)}$. The $W_d^{(i)}$ are pair-wise non-isogenous (since $\text{End}^{\! \circ \!}(W)$ is commutative), thus they are the isotypic components of $W$.

As above, one sees that $W_d^{(i)} = \Phi'^{(i)}_d(a)(W)$ and that $W_d^{(i)}$ is the reduced identity component of the kernel of $\Phi_d^{(i)}(a)$.

The component $W_d$ is simple if and only if $\Phi_d$ is irreducible over $\text{\rm End}^{\! \circ \!}(A)$, i.e.\ if and only if $\text{\rm End}^{\! \circ \!}(A) \otimes_{\mathbb{Q}} \mathbb{Q}(\zeta_d)$ is a field. 
If we fix an inclusion of $\text{\rm End}^{\! \circ \!}(A)$ into $\overline{\mathbb{Q}}$, this is the case if and only if $\text{\rm End}^{\! \circ \!}(A) \cap \mathbb{Q}(\zeta_d) = \mathbb{Q}$.

In particular, none of the $W_d$ splits further if $\text{\rm End}^{\! \circ \!}(A) = \mathbb{Q}$ as is the case if $A$ is an elliptic curve without complex multiplication (over $k$).

We proved:
\begin{theorem}
\label{W-factors}
Let $K|k$ be a cyclic field extension of degree $n$. Let $A$ be an abelian $k$-variety.

Let $W$ be the Weil restriction of $A_K$ with respect to $K|k$. For all $d|n$, $W$ contains canonically an abelian subvariety $W_d$ with ${W_d}_K \approx A_K^{\varphi(d)}$ (non-canonically), and $W$ is isogenous to the product of the $W_d$. Here, $W_1 = A$ itself.

Assume in addition that $A_K$ is simple, $\text{\rm End}^{\! \circ \!}(A_K)$ is commutative and all endomorphisms of $A_K$ can be defined over $k$. Then the isotypic components of $W$ are all simple, and its endomorphism rings are all commutative.  Fix an inclusion of $\text{\rm End}^{\! \circ \!}(A_K)$ into $\overline{\mathbb{Q}}$. Then for each $d$, $W_d$ is simple if and only if $\text{\rm End}^{\! \circ \!}(A) \cap \mathbb{Q}(\zeta_d) = \mathbb{Q}$.
\end{theorem}

Let $E$ be an ordinary elliptic $k$-curve with complex multiplication (over $k$). Fix an inclusion of $\text{End}^{\! \circ \!}(E)$ into $\overline{\mathbb{Q}}$.
As $\text{\rm End}^{\! \circ \!}(E)$ is a quadratic extension of $\mathbb{Q}$, $\text{\rm End}^{\! \circ \!}(E) \cap \mathbb{Q}(\zeta_d) \supsetneq \mathbb{Q}$ if and only if $\text{\rm End}^{\! \circ \!}(E) \subseteq \mathbb{Q}(\zeta_d)$. If this is the case then $\Phi_d$ splits into two polynomials of degree $\frac{1}{2} \varphi(d)$ over $\text{End}^{\! \circ \!}(E)$.

\begin{corollary}
\label{W-factors-elliptic-curve}
Under the assumptions of the theorem, let $E$ be an ordinary elliptic $k$-curve with $\text{\rm End}(E) = \text{\rm End}(E_K)$ (this condition is automatically satisfied over finite fields or if $n$ is odd).

Then for each $d$, $W_d$ is not simple if and only if $E$ has complex multiplication over $k$ and $\text{\rm End}^{\! \circ \!}(E) \subseteq \mathbb{Q}(\zeta_d)$. If this is the case, $W_d$ contains two simple non-isogenous abelian subvarieties of dimension $\frac{\varphi(d)}{2}$, and $W_d$ is isogenous to the product of these abelian subvarieties.

In partular, if $n$ is prime and $n \equiv 1$ mod $4$ or $n \equiv 3$ mod $4$ and $\sqrt{-n} \notin \text{\rm End}^{\! \circ \!}(E)$, we have an isogeny decomposition $W \sim E \times N$ where $N$ is simple.
\end{corollary}
\vspace{2ex}
\subsection*{Acknowledgements}
The results are part of the doctoral thesis of the first author. The thesis was supported by the Friedrich-Naumann-Stiftung, the Deutsche Forschungsgemeinschaft DFG and the Volkswagen Stiftung.

Special cases of some of the present results were obtained in the graduate thesis of the second author.

We would like to thank G.\ Frey for introducing us the the subject and for constant support during the preparation of our theses.

\vspace{1cm}
\noindent
Claus Diem: Institut f\"ur Experimentelle Mathematik, Universit\"at Duisburg-Essen, Essen, Germany. e-mail: diem@exp-math.uni-essen.de \\[2 ex]
Niko Naumann: Mathematisches Institut der Universit\"at M\"unster, M\"unster, Germany. e-mail: naumannn@uni-muenster.de

\end{document}